\documentclass[10pt,twoside]{article}
\usepackage{amsmath,amssymb}
\usepackage{amscd}

\newcommand{\ideal}{\mathfrak}

\newtheorem{thm}[subsection]{}

\DeclareMathOperator{\Spec}{Spec} 
\DeclareMathOperator{\Hom}{Hom} 
\DeclareMathOperator{\length}{length} 
\DeclareMathOperator{\Fr}{Fr} 
\DeclareMathOperator{\Frob}{Frob} 
\DeclareMathOperator{\Ver}{Ver} 
\DeclareMathOperator{\Ker}{Ker} 

\DeclareMathOperator{\rank}{rank}

\def\YEAR{\year}\newcount\VOL\VOL=\YEAR\advance\VOL by-1995
\def\firstpage{1}\def\lastpage{1000}
\def\received{}\def\revised{}
\def\communicated{}

\makeatletter
\def\magnification{\afterassignment\m@g\count@}
\def\m@g{\mag=\count@\hsize6.5truein\vsize8.9truein\dimen\footins8truein}
\makeatother

\font\eightrm=cmr8
\font\caps=cmcsc10                    
\font\Caps=cmcsc10 scaled \magstep1   

\makeatletter
\@specialpagefalse\headheight=8.5pt
\def\DocMath{}
\renewcommand{\@evenhead}{%
    \ifnum\thepage>\lastpage\rlap{\thepage}\hfill%
    \else\rlap{\thepage}\slshape\leftmark\hfill{\caps\SAuthor}\hfill\fi}%
\renewcommand{\@oddhead}{%
    \ifnum\thepage=\firstpage{\DocMath\hfill\llap{\thepage}}%
    \else{\slshape\rightmark}\hfill{\caps\STitle}\hfill\llap{\thepage}\fi}%
\makeatother

\def\TSkip{\bigskip}
\newbox\TheTitle{\obeylines\gdef\GetTitle #1
\ShortTitle  #2
\SubTitle    #3
\Author      #4
\ShortAuthor #5
\EndTitle
{\setbox\TheTitle=\vbox{\baselineskip=20pt\let\par=\cr\obeylines%
\halign{\centerline{\Caps##}\cr\noalign{\medskip}\cr#1\cr}}%
        \copy\TheTitle\TSkip\TSkip%
\def\next{#2}\ifx\next\empty\gdef\STitle{#1}\else\gdef\STitle{#2}\fi%
\def\next{#3}\ifx\next\empty%
    \else\setbox\TheTitle=\vbox{\baselineskip=20pt\let\par=\cr\obeylines%
    \halign{\centerline{\caps##} #3\cr}}\copy\TheTitle\TSkip\TSkip\fi%
\centerline{\caps #4}\TSkip\TSkip%
\def\next{#5}\ifx\next\empty\gdef\SAuthor{#4}\else\gdef\SAuthor{#5}\fi%
\ifx\received\empty\relax
    \else\centerline{\eightrm Received: \received}\fi%
\ifx\revised\empty\TSkip%
    \else\centerline{\eightrm Revised: \revised}\TSkip\fi%
\ifx\communicated\empty\relax
    \else\centerline{\eightrm Communicated by \communicated}\fi\TSkip\TSkip%
\catcode'015=5}}\def\Title{\obeylines\GetTitle}
\def\Abstract{\begingroup\narrower
    \parskip=\medskipamount\parindent=0pt{\caps Abstract. }}
\def\EndAbstract{\par\endgroup\TSkip}

\long\def\MSC#1\EndMSC{\def\arg{#1}\ifx\arg\empty\relax\else
     {\par\narrower\noindent%
     2000 Mathematics Subject Classification: #1\par}\fi}

\long\def\KEY#1\EndKEY{\def\arg{#1}\ifx\arg\empty\relax\else
        {\par\narrower\noindent Keywords and Phrases: #1\par}\fi\TSkip}

\newbox\TheAdd\def\Addresses{\vfill\copy\TheAdd\vfill
    \ifodd\number\lastpage\vfill\eject\phantom{.}\vfill\eject\fi}
{\obeylines\gdef\GetAddress #1
\Address #2 
\Address #3
\Address #4
\EndAddress
{\def\xs{4.3truecm}\parindent=0pt
\setbox0=\vtop{{\obeylines\hsize=\xs#1\par}}\def\next{#2}
\ifx\next\empty 
     \setbox\TheAdd=\hbox to\hsize{\hfill\copy0\hfill}
\else\setbox1=\vtop{{\obeylines\hsize=\xs#2\par}}\def\next{#3}
\ifx\next\empty 
     \setbox\TheAdd=\hbox to\hsize{\hfill\copy0\hfill\copy1\hfill}
\else\setbox2=\vtop{{\obeylines\hsize=\xs#3\par}}\def\next{#4}
\ifx\next\empty\ 
     \setbox\TheAdd=\vtop{\hbox to\hsize{\hfill\copy0\hfill\copy1\hfill}
                \vskip20pt\hbox to\hsize{\hfill\copy2\hfill}}
\else\setbox3=\vtop{{\obeylines\hsize=\xs#4\par}}
     \setbox\TheAdd=\vtop{\hbox to\hsize{\hfill\copy0\hfill\copy1\hfill}
                \vskip20pt\hbox to\hsize{\hfill\copy2\hfill\copy3\hfill}}
\fi\fi\fi\catcode'015=5}}\gdef\Address{\obeylines\GetAddress}

\begin{document}
\Title  Families of p-divisible groups with constant Newton polygon
\ShortTitle  p-divisible groups
\SubTitle   
\Author  Frans Oort and Thomas Zink
\ShortAuthor Oort and Zink
\EndTitle
\Abstract Let $X$ be a $p$-divisible group with constant Newton
polygon over a normal noetherian scheme $S$. We prove that there
exists an isogeny to $X \to Y$ such that $Y$ admits a slope
filtration. In case  $S$ is regular this was proved by N.Katz for dim
$S = 1$ and by T.Zink for dim $S \geq 1$. 
\EndAbstract
\MSC 14L05, 14F30
\EndMSC
\KEY 
\EndKEY
\Address  Frans Oort 
          Mathematisch Instituut
          P.O. Box. 80.010 
          NL - 3508 TA Utrecht 
          The Netherlands 
          email: oort@math.uu.nl 
\Address  Thomas Zink  
          Fakult\"at f\"ur Mathematik 
          Universit\"at Bielefeld
          Postfach 100131
          D-33501 Bielefeld
          Deutschland
          email: zink@mathematik.
          uni-bielefeld.de
\Address
\Address
\EndAddress

\section*{Introduction}
\noindent
\noindent
In this paper we work over base fields, and over base schemes over 
$\mathbb{F}_p$, i.e. we work entirely in characteristic $p$. 
We study $p$-divisible groups $X$ over a base scheme $S$ (and, colloquially, 
a $p$-divisible group over a base scheme of positive 
dimension will be called a ``family of $p$-divisible groups''), such
that the Newton polygon of a fiber $X_s$ is independent of $s \in
S$. We call $X$ a $p$-divisible group with constant Newton polygon.

\smallskip

A $p$-divisible group over a field has a slope filtration, see [Z1],
Corollary 13; for the definition of a slope filtration, see Definition 
\ref{Defslf}. Over a base of positive dimension a slope filtration can
only exist if the Newton polygon is constant. In Example \ref{exa1} we
show that even in this case there are $p$-divisible groups which do
not admit a slope filtration.

\smallskip

The main result of this paper has as a corollary that {\it for a 
$p$-divisible group with constant Newton polygon over a normal base
up to isogeny a slope filtration does exist}, see Corollary \ref{cor}.

\smallskip
 
We have access to this kind of questions by the definition
of a {\it completely
slope divisible} $p$-divisible group, see Definition \ref{Defcsld},
which implies a structure finer than a slope filtration. 
The main theorem of this paper, Theorem \ref{th}, says  that over 
a {\it normal} base this structure
on a $p$-divisible group exists up to isogeny. In [Z1], Theorem 7, 
this was shown to be
true over a regular base. In \ref{exa2} we
show that without the condition ``normal'' the conclusion of the
theorem does not hold.

\smallskip

Here is a motivation for this kind definition and of results:
\begin{itemize}
\item A $p$-divisible group over an
algebraically closed  field is isogenous with a $p$-divisible group which
can be defined over a finite field.
\item A $p$-divisible group over an algebraically closed  field is
completely slope divisible, if and only if
it is isomorphic with a direct sum of isoclinic $p$-divisible groups
which can be defined over a finite field, see \ref{charact2}.
\end{itemize} 
We see that a completely slope divisible $p$-divisible group comes 
``as close as possible'' to a constant one, 
in fact up to extensions of $p$-divisible groups annihilated by an
inseparable extension of the base, and up to monodromy.

\smallskip

{}From  Theorem  \ref{th} we deduce constancy results which generalize
results of Katz [K] and more recently of de Jong and Oort [JO]. In
particular we prove, Corollary \ref{7c} below:

\smallskip

{\it Let $R$ be a henselian local ring with residue field $k$. Let $h$
be a natural number. Then there exists a constant $c$ with the
following property. Let $X$ and $Y$ be isoclinic $p$-divisible 
groups over $S = \Spec R$ whose heights are smaller than $h$. 
Let $\psi : X_k \rightarrow Y_k$ be a homomorphism. Then $p^{c} \psi$
lifts to a homomorphism $X \rightarrow Y$.}

\section{Completely slope divisible $p$-divisible groups}
\noindent
In this section we present basic definitions and methods already 
used in the introduction.

\smallskip

Let $S$ be a scheme over $\mathbb{F}_p$.
Let $\Frob : S \rightarrow S$ be the absolute Frobenius morphism. 
For a scheme $G/S$ we write: 
\begin{displaymath}
   G^{(p)} \quad=\quad G \times_{S,\Frob} S.
\end{displaymath}
We denote by $\Fr = \Fr_G : G \rightarrow
G^{(p)}$ the Frobenius morphism relative to $S$. If $G$ is a finite
locally free commutative group scheme we write $\Ver = \Ver_G : G^{(p)}
\rightarrow G$ for the ``Verschiebung''.

\smallskip

Let $X$ be a $p$-divisible group over $S$. We denote by $X(n)$ the
kernel of the multiplication by $p^n : X \rightarrow X$. This is a
finite, locally free group scheme which has rank $p^{nh}$ if $X$ is of
height $h$.

\smallskip

Let $s = \Spec k $ the spectrum of a field of characteritic $p$. Let
$X$ be a $p$-divisible group over $s$. Let $\lambda \geq 0$ be a
rational number. We call $X$ isoclinic of slope $\lambda$, if there
exists integers $r \geq 0$,  $s > 0$ such that $\lambda = r/s$, and a
$p$-divisible group $Y$ over $S$, which is isogenous to $X$ such that
\begin{displaymath}
   p^{-r} \Fr^s : Y \rightarrow Y^{(p^s)}
\end{displaymath}
is an isomorphism.

\smallskip

A $p$-divisible group $X$ over $S$ is called isoclinic of slope
$\lambda $, if for each point $s \in S$ the group $X_s$
is isoclinic of slope $\lambda$. 

\vspace{-0.8cm}

\subsection*{}
\begin{thm}\label{Defslf} {\bf Definition.} Let $X/S$ be a
$p$-divisible group over a scheme $S$. A filtration  
$$
0 = X_0 \subset X_1 \subset X_2 \subset \ldots \subset X_m =X 
$$
consisting of $p$-divisible groups contained in $X$
is called a \emph{slope filtration}  of $X$ if there exists rational
numbers $\lambda_1, \ldots, \lambda_m$ satisfying $1 \geq \lambda_1 > 
\ldots > \lambda_m \geq 0$ such that every subquotient $X_i /
X_{i-1}$, $1<i \leq m$, is  isoclinic of slope $\lambda_i$.
\end{thm}

A $p$-divisible group $X$ over a field admits a slope filtration, see [Z1], 
Corollary 13. The slopes $\lambda_i$ and the heights of $X_i/X_{i-1} $
depend only on $X$. The height of $X_{i}/X_{i-1}$ is called the
multiplicity of $\lambda_i$. 
 
\smallskip
 
Over connected base scheme $S$ of positive dimension a slope
filtration of $X$ can only exist if the slopes of
$X_s$  and their multiplicities are independent of $s \in S$. In this
case we say that $X$ is a family of $p$-divisible groups with constant
Newton polygon.  Even if the Newton polygon is constant a slope
filtration in general does not exist, see Example \ref{exa1} below. 

\vspace{-0.8cm}

\subsection*{}
\begin{thm}\label{Defcsld} {\bf Definition.}   Let $s > 0$ and
$r_1, \ldots, r_m$ be integers such that $s \geq r_1 > r_2 > \ldots > 
r_m \geq 0$. A $p$-divisible group 
$Y$ over a scheme $S$ is said to be \emph{completely slope 
divisible} with respect to these integers if $Y$ has a filtration
 by $p$-divisible subgroups:
\begin{displaymath}
   0 = Y_0 \subset Y_1 \subset \ldots \subset Y_m = Y
\end{displaymath}
such that the following properties hold:
\begin{itemize}
\item The quasi-isogenies
\begin{displaymath}
      p^{-r_i}\Fr^s : Y_i \rightarrow Y_i^{(p^s)}
\end{displaymath} 
are isogenies for $i = 1, \ldots ,m$.
\item
The induced morphisms:
\begin{displaymath}
   p^{-r_i}\Fr^s : Y_i/Y_{i-1} \rightarrow (Y_i/Y_{i-1})^{(p^s)}
\end{displaymath}
are isomorphisms.
\end{itemize}
\end{thm}

Note that the last condition implies that $Y_i/Y_{i-1}$ is
isoclinic of slope $\lambda_i := r_i/s$. A filtration
described in this definition is a slope filtration in the sense of
the previous definition.

\bigskip
 
{\bf Remark.} Note that we do not require $s$ and $r_i$ to be 
relatively prime. If $Y$ is as in the definition, and $t \in \mathbb{Z}_{>0}$, 
it is also completely slope divisible with respect to $t{\cdot}s 
\geq t{\cdot}r_1 > t{\cdot}r_2 > \ldots > 
t{\cdot}r_m \geq 0$.  

We note that the filtration $Y_i$ of $Y$ is uniquely determined, if it exists.
Indeed, consider the isogeny $\Phi = p^{-r_m}\Fr^s: Y \rightarrow
Y^{(p^s)}$. Then $Y_m/Y_{m-1}$ is necessarily the $\Phi$-\'etale part
$Y^{\Phi}$ of  $Y$, see [Z1] respectively 1.6 below. This proves 
the uniqueness by induction.

\smallskip

We will say that a $p$-divisible group is completely slope divisible if it is 
completely slope divisible with respect to some set integers  and inequalities
$s \geq r_1 > r_2 > \ldots > r_m \geq 0$.

\bigskip

{\bf Remark.} A $p$-divisible group $Y$ over a field $K$ is completely slope 
divisible iff $Y \otimes_K L$ is completely slope 
divisible for some field $L \supset K$. - {\bf Proof.} The slope
filtration on $Y/K$ exists. We have $\Ker(p^{r}_Y) \subset
\Ker(\Fr^s_Y)$  iff  
$\Ker(p^{r}_{Y_L}) \subset \Ker(\Fr^s_{Y_L})$, and the same for equalities.
This proves that the conditions in the definition for completely slope 
divisibility hold over $K$ iff they hold over $L \supset K$. 

\vspace{-0.8cm}

\subsection*{}
\begin{thm}\label{perf} {\bf Proposition.} Let $Y$ be a completely slope 
divisible $p$-divisible group over a
perfect scheme $S$. Then $Y$ is isomorphic to a direct sum of 
isoclinic and completely slope divisible $p$-divisible  groups.
\end{thm}

\vspace{-0.2cm}

{\bf Proof.} With the notation of Definition \ref{Defcsld} we set
$\Phi = p^{-r_m}\Fr^{s}$. Let $Y(n) = \Spec \mathcal{A}(n)$ and let
$Y(n)^{\Phi} = \Spec \mathcal{L}(n)$ be the $\Phi$-\'etale part (see
Corollary \ref{GPhi} below). Then $\Phi$ induces a $\Frob^s$-linear
endomorphism $\Phi^{\ast}$ of $\mathcal{A}(n)$ and $\mathcal{L}(n)$
which is by definition bijective on $\mathcal{L}(n)$ and nilpotent on 
the quotient $\mathcal{A}(n)/ \mathcal{L}(n)$. One verifies (compare
[Z1], page 84) that there is a unique $\Phi^{\ast}$-equivariant section of
the inclusion $\mathcal{L}(n) \rightarrow \mathcal{A}(n)$.
This shows that $Y_m/Y_{m-1}$ is a direct factor of $Y = Y_m$. 
The result follows by induction on $m$.   \hfill $Q.E.D.$

\bigskip

Although not needed, we give a characterization of completely slope 
divisible $p$-divisible groups over a field. If $X$ is a $p$-divisible 
group over a field $K$, and $k \supset K$ is an algebraic closure, 
then $X$ is completely slope divisible if and only if $X_k$ is
completely slope divisible; hence it suffices to give a characterization over 
an algebraically closed field.
 
\smallskip

Convention: We will work with the covariant Dieudonn\'e module of a
$p$-divisible group over a perfect field ([Z2], [Me]). We write $V$,
respectively $F$ for Verschiebung, respectively Frobenius  on
Dieudonn\'e modules. Let $K$ be a perfect field, and let $W(K)$ be its
ring of Witt vectors. A Dieudonn\'e module $M$ over $K$ is the Dieudonn\'e
module of an isoclinic $p$-divisible group of slope $r/s$, iff there
exists a $W(K)$-submodule $M' \subset M$ such that $M/M'$ is
annihilated by a power of $p$, and such that $p^{-r}V^s (M') = M'$.

\smallskip

For later use we introduce the $p$-divisible group $G_{m,n}$ for
coprime positive integers $m$ and $n$. Its Dieudonn\'e module is
generated by one element, which is stable under  $F^m - V^n$. 
$G_{m,n}$  is isoclinic
of slope   $\lambda = m/(m+n)$ (in the terminology of this paper).
The height of $G_{m,n}$ is $h = m+n$, and this $p$-divisible group
is completely slope  divisible  with respect to $h=m+n \geq m$.
The group $G_{m,n}$ has dimension $m$, and its Serre dual 
has dimension $n$.

\smallskip

We have $G_{1,0} = \mu_{p^{\infty}}$, 
and $G_{0,1} = \underline{\mathbb{Q}_p/\mathbb{Z}_p}$. 

\vspace{-0.8cm}

\subsection*{}
\begin{thm}\label{charact1} {\bf Proposition.} Let  $k$ be an
algebraically closed field.  An isoclinic $p$-divisible group  $Y$
over  $k$ is completely slope 
divisible iff it can be defined over a finite field; i.e. 
iff there exists a $p$-divisible group  $Y'$ over some $\mathbb{F}_q$ 
and an isomorphism  $Y \cong Y' \otimes_{\mathbb{F}_q} k$.
\end{thm} 

\vspace{-0.2cm}

{\bf Proof.}  
Assume that $Y$ is slope divisible with respect to $s \geq r \geq 0$. 
Let $M$ be the covariant Dieudonn\'e module of $Y$. We set 
$\Phi = p^{-r} V^s$. By assumption this is a semilinear automorphism
of $M$. By a theorem of Dieudonn\'e (see 1.6 below) $M$ has a basis of 
$\Phi$-invariant vectors. Hence $M = M_0 \otimes_{W(\mathbb{F}_{p^s})}
W(k)$ where $M_0 \subset M$ is the subgroup of $\Phi$-invariant
vectors. Then $M_0$ is the Dieudonn\'e module of a $p$-divisible group
over $\mathbb{F}_{p^s}$ such that $Y \cong Y' \otimes_{\mathbb{F}_q}
k$.

\smallskip

Conversely assume that $Y$ is isoclinic over a finite field $k$ of
slope $r/s$. Let $M$ be the Dieudonn\'e module of $Y$. By definition
there is a finitely generated $W(k)$-submodule $M' \subset M \otimes
\mathbb{Q}$ such that $p^m M' \subset M \subset M'$ for some natural
number $m$, and such that $p^{-r}V^s (M') = M'$. Then $\Phi =
p^{-r}V^s$ is an automorphism of the finite set $M'/p^{m}M'$.  Hence
some power $\Phi^{t}$ acts trivially on this set. This implies that
$\Phi^t (M) = M$. We obtain that $p^{-rt}F^{st}$ induces an
automorphism of $Y$. Therefore $Y$ is completely slope divisible. 
\hfill $Q.E.D.$
\vspace{-0.8cm}

\subsection*{} 
\begin{thm}\label{charact2} {\bf Corollary.}  Let $Y$ be a
$p$-divisible group 
over an algebraically closed field $k$. This $p$-divisible group is
completely slope divisible iff $Y \cong \oplus Y_i$ such that every 
$Y_i$ is isoclinic, and can be defined over a finite field.
\end{thm}
{\bf Proof.} Indeed this follows from \ref{perf} and \ref{charact1}.
\hfill $Q.E.D.$

\subsection{The $\Phi$-\'etale part.} For further use, we recall a notion
explained and used in [Z1], Section 2. This method goes back to Hasse and 
Witt, see \cite{HW},
and to Dieudonn\'e, see \cite{D}, Proposition 5 on page 233. It can be 
formulated and proved for locally free sheaves, and it has a corollary for
finite flat group schemes.

\smallskip
 
Let $V$ be a finite dimensional vector space over a separably closed  
field $k$ of characteristic $p$. Let $f: V \to V$ be a $\Frob^{s}$-linear
endomorphism.  The set $C_V = \{x \in V \mid
f(x) = x\}$ is a vectorspace  over $\mathbb{F}_{p^s}$. Then $V^{f} = C_V
\otimes_{\mathbb{F}_{p^s}} k$ is a subspace of $V$. The endomorphism
$f$ acts as a $\Frob^s$-linear automorphism on $V^{f}$ and acts
nilpotently on the quotient $V/V^{f}$. This follows essentially form 
Dieudonn\'e loc.cit.. Moreover if $k$ is any field
of characteristic $p$ we have still unique exact sequence of $k$-vector
spaces
\begin{displaymath}
   0 \rightarrow V^{f} \rightarrow V \rightarrow V/V^{f} \rightarrow 0
\end{displaymath}
such that $f$ acts as a a $\Frob^s$-linear automorphism on $V^{f}$ and acts
nilpotently on the quotient $V/V^{f}$.

\smallskip
 
This can be applied in the following situation:  
Let $S$ be a scheme over $\mathbb{F}_p$. 
Let $G$ be a locally free group scheme over $S$ endowed with a
homomorphism 
$$\Phi: G \to G^{(p^s)}.$$ 
 
\smallskip
 
In case $S = \Spec(K)$, where $K$ is a field we consider the affine
algebra $A$ of $G$. The $\Phi$ induces a $\Frob^s$-linear endomorphism $f: A
\rightarrow A$. The vector subspace $A^f$ inherits the structure of a
bigebra. We obtain a finite group scheme $G^{\Phi} = \Spec A^{f}$,
which is called the $\Phi$-\'etale part of $G$. Moreover we have an
exact sequence of group schemes:
\begin{displaymath}
 0 \to G^{\Phi-{\rm nil}} \longrightarrow  G  \longrightarrow G^{\Phi} \to 0.  
\end{displaymath}
The morphism $\Phi$ induces an isomorphism $G^{\Phi} \rightarrow
(G^{\Phi})^{(p^s)}$, and acts nilpotently on the kernel $G^{\Phi-{\rm nil}}$.

\smallskip
 
Let now $S$ be an arbitrary scheme over $\mathbb{F}_p$. Then we can
expect a $\Phi$-\'etale part only in the case where the rank of
$G_x^{\Phi}$ is independent of $x \in S$:

\vspace{-0.8cm}

\subsection*{}
\begin{thm}\label{GPhi} {\bf Corollary.}  Let $G \to S$ be a finite,
locally free  group scheme; let $\Phi: G \to G^{(p^t)}$ be a
homomorphism.  Assume that  the function 
$$S \to \mathbb{Z}, \quad\mbox{\it defined by}\quad x \mapsto
\rank((G_x)^{\Phi_x})$$ 
is constant. Then there exists an exact sequence
$$0 \to G^{\Phi-{\rm nil}} \longrightarrow  G  \longrightarrow G^{\Phi} \to 0$$
such that $\Phi$ is nilpotent on $G^{\Phi-{\rm nil}}$ and an isomorphism on
the $\Phi$-\'etale part $G^{\Phi}$.
\end{thm}

The prove is based on another proposition which we use in section 3.
Let $\mathcal{M}$ be a finitely generated, locally 
free $\mathcal{O}_S$-module. Let  $t \in \mathbb{Z}_{>0}$
$$f: \mathcal{O}_S \otimes_{\Frob_S^t, S} \mathcal{M} =
\mathcal{M}^{(p^t)} \rightarrow \mathcal{M}$$ 
be a morphism of $\mathcal{O}_S$-modules.  To every morphism $T \to S$ we
associate 
$$C_{\mathcal{M}}(T) = \{x \in \Gamma(T,\mathcal{M}_T) \mid f(1
\otimes x) = x\}.$$

\vspace{-0.8cm}

\subsection*{}
\begin{thm}\label{etfin} {\bf Proposition} ({\rm see \cite{Z1}, Proposition 3}).
The functor $C_{\mathcal{M}}$ is 
represented by a scheme that is \'etale and affine over $S$. Suppose $S$ 
to be connected; the
scheme $C_{\mathcal{M}}$ is finite over $S$ iff for each geometric point $\eta
\rightarrow S$ the cardinality of $C_{\mathcal{M}}(\eta)$ is the same. 
\end{thm}

Let $X$ be a $p$-divisible group over a field $K$. Suppose $\Phi : X
\rightarrow X^{(p^t)} $ is a homomorphism. Then the $\Phi$-\'etale
part $X^{\Phi}$ is the inductive limit of $X(n)^{\Phi}$. This is a
$p$-divisible group. 

\vspace{-0.8cm}

\subsection*{}
\begin{thm}\label{nil} {\bf Corollary.}  Let $X$ be a
$p$-divisible group over $S$. Assume that for each geometric point
$\eta \rightarrow S$ the height of the $\Phi$-\'etale part of
$X_{\eta}$ is the same. Then a $p$-divisible group $X^{\Phi}$ exists
and commutes with arbitrary base change. There is an exact sequence of
$p$-divisible groups: 
\begin{displaymath}
    0 \rightarrow X^{\Phi-{\rm nil}} \rightarrow X \rightarrow X^{\Phi}
   \rightarrow 0. 
\end{displaymath}
\end{thm}

The following proposition can be deduced from proposition
\ref{etfin}. 

\vspace{-0.8cm}

\subsection*{}
\begin{thm}\label{Dieurel} {\bf Corollary.}  Assume that $G \to
S$ is a finite, 
locally free group scheme over a connected base scheme $S$.  
Let $\Phi: G \stackrel{\sim}{\longrightarrow} G^{(q)}$,  $q = p^s$ be 
an isomorphism. Then there exists a finite \'etale morphism $T \to S$,
and a morphism $T \to \Spec(\mathbb{F}_{q})$, 
such that $G_T$ is obtained by base change from a finite group scheme
$H$ over $\mathbb{F}_{q}$:   
$$H \otimes_{\Spec \mathbb{F}_{q}} T  \stackrel{\sim}{\longrightarrow}
G_T.$$  
Moreover $\Phi$ is induced from the identity on $H$. 
\end{thm}

{\bf Remark.} If $S$ is a scheme over $\bar{\mathbb{F}}_p$ the
Corollary says in 
particular that $G_T$ is obtained by base change from a finite group
scheme over $\bar{\mathbb{F}}_p$. In this case we call $G_T$ constant
(compare \cite{K}, (2.7)). This should not be confused with the
\'etale group scheme associated to a finite abelian group $A$. 
We will discuss ``constant'' $p$-divisible groups, see Section 3 below.

\section{The main result: slope filtrations}
\noindent
In this section we show:
\vspace{-0.8cm}

\subsection*{} 
\begin{thm}\label{th} {\bf Theorem.}  Let $h$ be a natural
number. Then there exists a natural number $N(h)$ with the following
property. Let $S$ be an integral, normal  noetherian scheme.
Let $X$ be a
p-divisible group over $S$ of height $h$ with constant Newton
polygon. Then there is  a completely slope divisible $p$-divisible
group $Y$ over $S$, and an isogeny:
\begin{displaymath}
 \varphi: X \rightarrow Y \quad\mbox{\it over}\quad S 
 \quad\mbox{\it with}\quad \deg(\varphi) \leq N(h).
\end{displaymath}
\end{thm}

\noindent
In Example \ref{exa2} we see that the condition ``$S$ is normal'' is essential.
By this theorem we see that a slope filtration exists up to isogeny:
\vspace{-0.8cm}

\subsection*{}
\begin{thm}\label{cor} {\bf Corollary.}  Let $X$ be a $p$-divisible group 
with constant Newton polygon over an integral, normal  noetherian
scheme $S$. There exists an isogeny 
$\varphi: X \to Y$, such that $Y$ over $S$ admits a slope filtration.
\end{thm}
\vspace{-0,5cm} 
\hfill $Q.E.D.$
 
\vspace{-0.8cm}

\subsection*{}
\begin{thm}\label{3p} {\bf Proposition.} Let $S$ be an integral scheme 
with function field $K = \kappa (S)$. Let $X$ 
be a $p$-divisible group over $S$ with constant Newton polygon, such that 
$X_K$ is completely slope divisible with respect to the integers 
$s \geq r_1 > r_2 > \ldots > r_m \geq 0$. 
Then $X$ is completely slope divisible with respect to the same
integers.
\end{thm}
{\bf Proof.}  The quasi-isogeny $\Phi = p^{-r_m}\Fr^s : X \rightarrow
X^{(p^s)}$ is an isogeny, because this is true over the general
point. Over any geometric point $\eta \rightarrow S$ the
$\Phi$-\'etale part of $X_{\eta}$ has the same height by constancy of
the Newton polygon. Hence the $\Phi$-\'etale part of $X$ exists by
Corollary \ref{nil}. We obtain an exact sequence:
\begin{displaymath}
0 \rightarrow X^{\Phi-{\rm nil}} \rightarrow X \rightarrow X^{\Phi}
\rightarrow 0.
\end{displaymath} 
Assuming  an induction hypothesis on  $X^{\Phi-{\rm nil}}$ gives the result.
\hfill $Q.E.D.$

\bigskip

A basic tool in the following proofs is the moduli scheme of isogenies
of degree $d$ of a $p$-divisible group (compare [RZ] 2.22): Let $X$ be
a $p$-divisible group over a scheme $S$, and let $d$ be a natural
number. Then we define the following functor $\mathcal{M}$ on the
category of $S$-schemes $T$. A point of $\mathcal{M}(T)$ consists of a
$p$-divisible group $Z$ over $T$ and an isogeny $\alpha : X_T
\rightarrow Z$ of 
degree $d$ up to isomorphim. \emph{The functor $\mathcal{M}$ is
representable by a projective scheme over $S$.}  Indeed, to each
finite, locally free 
subgroup scheme $G \subset X_T$ there is a unique isogeny $\alpha $
with kernel $G$. Let $n$ be a natural number such that $p^n \geq
d$. Then $G$ is a finite, locally free subgroup scheme on $X(n)_T$. We
set $X(n) = \Spec_S \mathcal{A}$. The affine algebra of $G$ is a
quotient of the locally free sheaf $\mathcal{A}_T$. Hence we obtain a
point of the Grassmannian of $\mathcal{A}$. This proves that
$\mathcal{M}$ is representable as a closed subscheme of this
Grassmannian.

\vspace{-0.8cm}

\subsection*{}
\begin{thm}\label{4l} {\bf Lemma.}  For every $h \in \mathbb{Z}_{>0}$
there  
exists a number $N(h) \in  \mathbb{Z}$ with the following property. Let
$S$ be an integral noetherian scheme. Let $X$ be a $p$-divisible group
of height $h$ over $S$ with constant Newton polygon.
There is a non-empty open subset $U \subset S$, and a projective morphism
$\pi : S^{\sim} \rightarrow S$ of integral schemes which induces an 
isomorphism $\pi: \pi^{-1}(U) \rightarrow U$ such that there exist a 
completely slope divisible $p$-divisible group $Y$ over $S^{\sim}$, and 
an isogeny $ X_{S^{\sim}} \rightarrow Y$, whose degree is bounded by
$N(h)$.
\end{thm}
{\bf Proof.}   Let $K$ be the function field of $S$. 
We know by [Z1], Prop. 12, that there is a completely slope divisible 
$p$-divisible group $Y^{0}$ over $K$, and an isogeny 
$\beta^{0} : X_K \rightarrow Y^{0}$, whose degree is bounded by 
a constant which depends only on the height of $X$. The kernel of this 
isogeny is a finite group scheme $G^{0} \subset X_K(n)$, for some $n$. 
Let $\bar{G}$ be the scheme-theoretic image of $G^{0}$ in $X(n)$, 
see EGA, I.9.5.3. 
Then $\bar{G}$ is flat over some nonempty open set $U \subset S$, and 
inherits there the structure of a finite, locally free group scheme $G 
\subset X(n)_U$. We form the $p$-divisible group $Z = X_U/G$. By
construction there are integers $s \geq r_m$, such that 
\begin{displaymath}
p^{-r_m}\Fr^s : Z_K \rightarrow Z_K^{(p^s)} 
\end{displaymath}
is an isogeny, and $r_m/s$ is a smallest slope in the Newton polygon
of $X$.  

Therefore $\Phi = p^{-r_m}\Fr^s : Z \rightarrow Z^{(p^s)}$ 
is an isogeny too. As in the proof of the last proposition 
the constancy of the Newton polygon implies that the $\Phi$-\'etale
part $Z^{\Phi}$ exists. We obtain an 
exact sequence of $p$-divisible groups on $U$:
\begin{displaymath}
0 \rightarrow Z^{\Phi-{\rm nil}} \rightarrow Z \rightarrow Z^{\Phi}
\rightarrow 0.  
\end{displaymath} 
By induction we find a non-empty open subset $V \subset U$ and a 
completely slope divisible $p$-divisible group $Y_{m-1}$ which is 
isogenous to $Z^{\Phi-{\rm nil}}_V$. Taking the push-out of the last exact 
sequence by the isogeny $Z^{\Phi-{\rm nil}}_V \rightarrow Y_{m-1}$ we find
a completely slope divisible $p$-divisible group $Y$ over $V$ which is
isogenous to $X_V$. 
 
\smallskip

Let $d$ be the degree of the isogeny $\rho : X_V \rightarrow Y$. 
We consider the 
moduli scheme $\mathcal{M}$ of isogenies of degree $d$ of $X$ defined
above. 
The isogeny $X_V \rightarrow Y$ defines an $S$-morphism $V \rightarrow 
\mathcal{M}$. The scheme-theoretic image $S^{\sim}$ of $V$   
is an integral scheme, which is projective over $S$. Moreover
the morphism $\pi : S^{\sim} \rightarrow S$ induces an isomorphism 
$\pi^{-1} (V) \rightarrow V$. The closed immersion $S^{\sim} \rightarrow 
\mathcal{M}$ corresponds to an isogeny $ \rho^{\sim}: 
X_{S^{\sim}} \rightarrow Y^{\sim}$ to a $p$-divisible group 
$Y^{\sim}$ on $S^{\sim}$. Moreover the restriction of $\rho^{\sim}$ to  
$V$ is $\rho$. Since $Y^{\sim}$ has constant Newton polygon, and since 
$Y^{\sim}$ is completely slope divisible in the generic point of $S^{\sim}$ 
it is completely slope divisible by Proposition \ref{3p}.       
 \hfill $Q.E.D.$
\vspace{-0.8cm}

\subsection*{}
\begin{thm}\label{6l} {\bf Lemma.}   Let $k$ be an 
algebraically closed field of characteristic $p$. Let
$s \geq r_1 > r_2 > \ldots > r_m \geq 0$ and $d > 0$ be integers.
Let $X$ be a $p$-divisible group over $k$. Then there are up to 
isomorphism
only finitely many isogenies $X \rightarrow Z$ of degree $d$ to a 
$p$-divisible group $Z$, which is completely slope divisible with respect 
to $s \geq r_1 > r_2 > \ldots > r_m \geq 0$.
\end{thm}
{\bf Proof.} It suffices to show this in case also $X$ is 
completely slope divisible with 
respect to $s \geq r_1 > r_2 > \ldots > r_m \geq 0$. Then $X$ and $Z$
are a direct product of isoclinic slope divisible groups. Therefore
we assume that we are in the isoclinic case $m = 1$.

\smallskip

In this case we consider the contravariant Dieudonn\'e modules $M$ of $X$, 
and $N$ of $Z$. Let $\sigma $ be the Frobenius on $W(k)$. 
Then $N \subset M$ is a submodule such that $\length M/N = \log_p d$. By 
assumption $\Phi = p^{-r_1}F^s$ induces a $\sigma^s$-linear automorphism of 
$M$ respectively $N$. Let $C_N$ respectively $C_M$ be the invariants
of $\Phi$ acting of $N$ respectively $M$. Hence $C_N$ is a
$W(\mathbb{F}_{p^s})$-submodule of $N$ such that $W(k)
\otimes_{W(\mathbb{F}_{p^s})} C_N = N$ (e.g. \cite{Z2} 6.26). The same
holds for $M$.  
We see that $C_N$ is a $W(\mathbb{F}_{p^s})$-submodule of $C_M$, such
that $\length C_M/C_N = \log_p d$. Since there are only finitely many
such submodules, the assertion follows.  
\hfill $Q.E.D.$
\vspace{-0.8cm}

\subsection*{}
\begin{thm}\label{l5} {\bf Lemma.}  Let
$f: T \rightarrow S$ be a proper morphism of  schemes
such that $f_{\ast} \mathcal{O}_T = \mathcal{O}_S$. Let $g : T
\rightarrow M$ be a morphism of schemes. 
We assume that for any point $\xi \in S$ the set-theoretic
image of the fiber $T_{\xi}$ by $g$ is a single point in $M$. 
Then there is a unique morphism $h : S \rightarrow M$ such that $hf =
g$.
\end{thm}
{\bf Proof.}  For $\xi \in S$ we set $h(\xi) = f(T_{\xi})$. This defines a 
set-theoretic map $h : S \rightarrow M$. If $U \subset M$ is an open
neighborhood of  
$h(\xi)$ then $g^{-1}(U)$ is an open neighborhood of $T_{\xi}$. Since $f$
is closed we find an open neighbourhood $V$ of $\xi$ with $f^{-1}(V) \subset 
g^{-1}(U)$. Hence $h(V) \subset U$; we see that $h$ is continuous. Then 
$h_{\ast} \mathcal{O}_S =  h_{\ast} f_{\ast} \mathcal{O}_T = 
g_{\ast} \mathcal{O}_T$. We obtain a morphism of ringed spaces 
$\mathcal{O}_M \rightarrow g_{\ast} \mathcal{O}_T = h_{\ast} \mathcal{O}_S$.
\hfill  $Q.E.D.$

\vspace{5mm} \noindent
Theorem \ref{th} follows from the following technical variant which is
useful if we do not know that the normalization is finite. We will
need that later on. 
\vspace{-0.8cm}

\subsection*{}
\begin{thm}\label{thmneu} {\bf Proposition.}  Let $h$ be a natural
number. Then there exists a natural number $N(h)$ with the following
property. Let $S$ be an integral noetherian scheme. Let $X$ be a
p-divisible group over $S$ of height $h$ with constant Newton
polygon. Then there is a finite birational morphism 
$T \rightarrow S$, a completely slope divisible $p$-divisible
group $Y$ over $T$, and an isogeny:
\begin{displaymath}
  X_{T} \rightarrow Y
\end{displaymath}
over $T$ whose degree is smaller than $N(h)$. 
\end{thm}
{\bf Proof.}   Consider the proper birational map 
$\pi : S^{\sim} \rightarrow S$, and the isogeny 
$\rho :X_{S^{\sim}} \rightarrow Y$ given by 
Lemma \ref{4l}. Take the Stein factorization $S^{\sim} \rightarrow 
T \rightarrow S$. It is enough to find over $T$ an isogeny to a completely 
slope divisible $p$-divisible group. Therefore we  assume $S = T$, i.e.
$\pi_{\ast} \mathcal{O}_{S^{\sim}} = \mathcal{O}_S$.

\smallskip

Let  $\mathcal{M} \to S$ be the moduli scheme of isogenies of 
$X$ of degree $d = \mbox{degree} \; \rho$. We will show that the
$S$-morphism $g: S^{\sim} \rightarrow \mathcal{M}$ defined by $Y$
factors through $S \rightarrow \mathcal{M}$.

\smallskip

Let $\xi \in S$. We write 
$S^{\sim}_{\xi} = S^{\sim} \times_{S} \Spec(\kappa(\xi))$. 
By Lemma \ref{l5} it is suffices to show that the 
set-theoretic image of $S^{\sim}_{\xi}$ by 
$g$ is a single point of $\mathcal{M}$. Clearly $\mathcal{M}_{\xi}$
classifies isogenies starting at  $X_{\xi}$ of degree $d$. 
Over the algebraic closure $\bar{\xi}$,  by Lemma \ref{6l}, 
there are only finitely many
isogenies of $X_{\bar{\xi}} \rightarrow Z$ of degree $d$ to a
completely slope divisible group (for fixed $s$ and $r_i$).
This shows that the image of
$S^{\sim}_{\xi} (\bar{\xi}) \rightarrow \mathcal{M}_{\xi} (\bar{\xi})$ is
finite. Since $S^{\sim}_{\xi}$ is  connected, see [EGA], III$^1$.4.3.1, the
image of $S^{\sim}_{\xi}$ is a single point of $\mathcal{M}_{\xi}$.

\smallskip

Hence we have the desired factorization $S \rightarrow \mathcal{M}$. 
It defines an isogeny $X \rightarrow Z$ over $S$. Finally 
$Z$ is completely 
slope divisible since it is completely slope divisible in the general point
of $S$, and because its Newton polygon is constant. 

\hfill $Q.E.D.$ \ref{thmneu} \& \ref{th}

\section{Constancy results}
\noindent
Let $T$ be a scheme over $\bar{\mathbb{F}}_p$. We study the question if a
$p$-divisible group $X$ over $T$ is constant up to isogeny, i.e. there 
exists a $p$-divisible group $Y$ over $\bar{\mathbb{F}}_p$ such that $X$ is 
isogenous to $Y \times_{\bar{\mathbb{F}}_p} T$.
\vspace{-0.8cm}

\subsection*{}
\begin{thm}\label{5c2} {\bf Proposition.} Let $S$ be a noetherian
integral normal scheme over $\bar{\mathbb{F}}_p$. Let 
$K$ be the function field of $S$ and let ${\bar K}$ be an algebraic
closure of $K$. We denote by $L \subset {\bar K}$ the maximal 
unramified extension of $K$ with respect to $S$. Let $T$ be
the normalization of $S$ in $L$. 

Let $X$ be an isoclinic $p$-divisible group over $S$. Then there is a
$p$-divisible group $X_0$ over $\bar{\mathbb{F}}_p$ and an isogeny
$  X \times_{S} T \rightarrow  X_0 \times_{\Spec \bar{\mathbb{F}}_p} T
$ such that  the degree of this isogeny is smaller than an integer
which depends only on the height of $X$.  
\end{thm}
{\bf Proof:} We use Theorem \ref{th}: there exists an 
isogeny $\varphi: X \to Y$, where $Y$ over $S$ is completely slope divisible. 
There are natural numbers $r$ and $s$, such that 
\begin{equation*}
\Phi = p^{-r} \Fr^s : Y \rightarrow Y^{(p^s)}
\end{equation*}
is an isomorphism. Applying Corollary \ref{Dieurel} to $Y(n)$ and $\Phi$ we
obtain finite group schemes $X_0(n)$ over $\mathbb{F}_{p^s}$ and isomorphisms
\begin{displaymath}
   Y(n)_T \cong X_0(n) \times_{\mathbb{F}_{p^s}} T
\end{displaymath}
The inductive limit of the group schemes $X_0(n)$ is a $p$-divisible
group $X_0$ over $\mathbb{F}_{p^s}$. It is isogenous to $X$ over $T$.
  \hfill $Q.E.D.$
\vspace{-0.8cm}

\subsection*{} 
\begin{thm}{\bf Corollary.} Let $S$ and $T$ be as in the
proposition. Let $T^{\rm perf} \rightarrow T$ be the perfect hull of 
$T$. 
Let $X$ be a $p$-divisible group over $S$ with constant Newton
polygon. Then there is a
$p$-divisible group $X_0$ over $\bar{\mathbb{F}}_p$ and an isogeny
$X_0 \times_{\Spec \bar{\mathbb{F}}_p} T^{\rm perf} \rightarrow X
\times_{S} T^{\rm perf}$, whose degree is smaller than an integer which
depends only on the height of $X$.  
\end{thm}
{\bf Proof.} This follows using Proposition \ref{perf}. 
\hfill $Q.E.D.$

\medskip

Finally we prove constancy results without the normality condition. 
\vspace{-0.8cm}

\subsection*{}
\begin{thm}\label{6p} {\bf Proposition.} 
Let $R$ be a strictly henselian reduced local ring over
$\bar{\mathbb{F}}_p$. Let $X$ be an isoclinic $p$-divisible group over
$S = \Spec R$. Then there is a
$p$-divisible group $X_0$ over $\bar{\mathbb{F}}_p$ and an isogeny
$X_0 \times_{\Spec \bar{\mathbb{F}}_p} S \rightarrow X$, whose degree
is smaller than an integer which depends only on the height of $X$.
\end{thm}
\vspace{-1.2cm}

\subsection*{}
\begin{thm}\label{7c} {\bf Corollary.} 
To each natural number $h$ there is a natural number $c$ with the
following property:
Let $R$ be a henselian reduced local ring over $\mathbb{F}_p$ with
residue field $k$. Let $X$ and $Y$ be isoclinic $p$-divisible 
groups over $S = \Spec R$ whose heights are smaller than $h$. 
Let $\psi : X_k \rightarrow Y_k$ be a homomorphism. Then $p^c \psi$
lifts to a homomorphism $X \rightarrow Y$.
\end{thm}
 
\noindent
A proof of the proposition, and of the 
corollary will be given later. 

\medskip

{\bf Remark.} In case the $R$ considered in the previous proposition,
or in the previous corollary, is not reduced, but satisfies all other
properties, the conclusions still hold, except that the integer
bounding the degree of the  
isogeny, respectively the integer $c$, depend on $h$ and on $R$. 

\medskip

If $R$ is strictly henselian the corollary follows from Proposition
\ref{6p}. Indeed, assume that $X$ and $Y$ are isogenous to constant
$p$-divisible groups $X_0$ and $Y_0$  by isogenies which
are bounded by a constant which depends only on $h$. The corollary
follows because:
\begin{displaymath}
   \Hom ((X_0)_k, (Y_0)_k) = \Hom ((X_0)_R, (Y_0)_R).
\end{displaymath}
Conversely the corollary implies the proposition since by Proposition 
\ref{5c2} over a separably closed field an isoclinic $p$-divisible
group is isogenous to a constant $p$-divisible group.  

\medskip

{\bf Remark.}
Assume Corollary \ref{7c}.
An isoclinic slope divisible $p$-divisible group $Y$ over $k$ can be
lifted to an isoclinic slope divisible $p$-divisible group over
$R$. Indeed the \'etale schemes associated by \ref{etfin} 
to the affine algebra of $Y(n)$ and the isomorphism $p^{-r}\Fr^{s}$
lift to $R$. Hence the categories of isoclinic $p$-divisible groups up
to isogeny over $R$ respectively $k$ are equivalent. 

\vspace{-0.8cm}

\subsection*{}
\begin{thm}\label{7l} {\bf Lemma.} 
Consider a commutative diagram of rings over $\mathbb{F}_p$:
\begin{equation*}
\begin{array}{ccc}
R & \rightarrow & A\\
\downarrow &  & \downarrow\\
R_0 & \rightarrow & A_0.
\end{array}
\end{equation*}
Assume that $R \rightarrow R_0$ is a surjection with nilpotent kernel
$\ideal{a}$, and that $A \rightarrow A_0$ is a surjection with
nilpotent kernel $\ideal{b}$. Moreover let $R \rightarrow A$ be a
monomorphism. 

Let $X$ and $Y$ be $p$-divisible groups over $R$. Let
$\varphi_0 :X_{R_0} \rightarrow Y_{R_0}$ be a morphism of the $p$-divisible
groups obtained by base change. Applying base change with respect to
$R_0 \rightarrow A_0$ we obtain a morphism $\psi_0 : X_{A_0}
\rightarrow Y_{A_0}$. 

If $\psi_0$ lifts to a morphism $\psi : X_A \rightarrow Y_A$, then
$\varphi_0$ lifts to a morphism $\varphi : X \rightarrow Y$.
\end{thm}
{\bf Proof.} By rigidity, liftings of homomorphisms of $p$-divisible
groups are unique. Therefore we may replace $R_0$ by its image in
$A_0$ and assume that $R_0 \rightarrow A_0$ is injective. Then we
obtain $\ideal{a} = \ideal{b} \cap R$.

\smallskip

Let $n$ be a natural number such that $\ideal{b}^n = 0$. We
argue by induction on $n$. If $n = 0$, we have $\ideal{b} = 0$ and 
therefore $\ideal{a} = 0$. In this case there is nothing to prove. 
If $n > 0$ we consider the commutative diagram:
\begin{equation*}
\begin{array}{ccc}
R & \rightarrow & A\\
\downarrow &  & \downarrow\\
R/(\ideal{b}^{n-1} \cap R) & \rightarrow & A/\ideal{b}^{n-1}\\
\downarrow &  & \downarrow\\
R_0 & \rightarrow & A_0.
\end{array}
\end{equation*}
We apply the induction hypothesis to the lower square. Hence it is 
enough to show the lemma for the upper square. We assume therefore
without loss of generality that $\ideal{a}^2 = 0, \; \ideal{b}^2 = 0$.

\smallskip

Let $D_X$ and $D_Y$ be the crystals associated to $X$ and $Y$ by
Messing [Me]. The values $D_X (R)$ respectively $D_Y(R)$ are finitely
generated projective $R$-modules which are endowed with the Hodge
filtration $Fil_X \subset D_X (R)$ respectively $Fil_Y \subset
D_Y(R)$. We put on $\ideal{a}$ respectively $\ideal{b}$ the trivial
divided power structure. Then $\varphi_0 $ induces a map $D_X(R)
\rightarrow D_Y(R)$. By the criterion of Grothendieck and Messing
$\varphi_0 $ lifts to a homomorphism over $R$, iff $D(\varphi_0)
(Fil_X ) \subset Fil_Y$.

\smallskip

Since the construction of the crystal commutes with base change, see [Me],
Chapt. IV, 2.4.4, we have canonical isomorphisms:
\begin{displaymath}
   \begin{array}{cc}
     D_{X_A}(A) = A \otimes_R D_X (R), \quad  & Fil_{X_A} = A \otimes_R
     Fil_X,\\
     D_{Y_A}(A) = A \otimes_R D_Y (R), & Fil_{Y_A} = A \otimes_R
     Fil_Y.\\
   \end{array}
\end{displaymath}
Since $\psi_0$ lifts we have $id_A \otimes D(\varphi) (A \otimes_R Fil_X)
\subset A \otimes_R Fil_Y$. Since $R \rightarrow A$ is injective this
implies $D(\varphi_0) (Fil_X) \subset Fil_Y$.  
\hfill  $Q.E.D.$

\medskip

{\bf Proof} of Proposition \ref{6p}. We begin with the case where $R$
is an integral domain.  
By Proposition \ref{thmneu} there is a
finite ring extension $R \rightarrow A$ such that $A$ is contained in 
the quotient field of $R$, and such that there is an isogeny $X_A
\rightarrow  Y$ to a completely slope divisible $p$-divisible group $Y$ over
$A$. The degree of this isogeny is smaller than a constant
which depends only on the height of $X$. Since $A$ is a product of
local rings we may assume without loss of generality that $A$ is
local. The ring $A$ is a strictly henselian local ring, see [EGA] IV
18.5.10, and has
therefore no non-trivial finite \'etale coverings. The argument of the proof of
Proposition \ref{5c2} shows that $Y$ is obtained by base change from a 
$p$-divisible group $X_0$ over $\bar{\mathbb{F}}_p$. Therefore we find
an isogeny 
\begin{displaymath}
   \varphi : (X_0)_A \rightarrow X_A
\end{displaymath}
Let us denote by $k$ the common residue field of $A$ and $R$. Then
$\varphi$ induces an isogeny $\varphi: (X_0)_k \rightarrow X_k$. 
The last lemma shows that $\varphi_0$ lifts to an isogeny
$\hat{\varphi} : (X_0)_{\hat{R}} \rightarrow X_{\hat{R}}$ over the
completion $\hat{R}$ of $R$. 

\smallskip
 
We apply the following  fact: 

\noindent 
{\bf Claim.} Consider a fiber product of rings: 
\begin{equation*}
\begin{array}{ccc}
R & \rightarrow & A_1\\
\downarrow & & \downarrow\\
A_2 & \rightarrow & B.
\end{array}   
\end{equation*}
Let $X$ and $Y$ be $p$-divisible groups over $R$. Let $\psi_i :
X_{A_i} \rightarrow Y_{A_i}$ for $i = 1, 2$ be two homomorphisms of
$p$-divisible groups which agree over $B$. Then there is a
unique homomorphism $\psi : X \rightarrow  Y$ which induces $\psi_1$
and $\psi_2$.

\smallskip

In our concrete situation we consider the diagram:
\begin{equation*}
\begin{array}{ccc}
R & \rightarrow & A\\
\downarrow & & \downarrow\\
\hat{R} & \rightarrow & \hat{A} = \hat{R} \otimes_R A.
\end{array}   
\end{equation*}
The morphisms $\hat{\varphi}$ and $\varphi$ agree 
over $\hat{A}$ because they agree over the residue field $k$. 
This proves the case of an integral domain $R$.

\smallskip

In particular we have shown the Corollary \ref{7c} in the case where
$R$ is a strictly henselian integral domain. To show the corollary in
the reduced case 
we consider the minimal prime ideals $\ideal{p}_1, \ldots \ideal{p}_s$
of $R$. Let $\psi : X_k \rightarrow  Y_k$ be a homomorphism. Then $p^c
\psi$ lifts to a homomorphims over each of the rings $R/\ideal{p}_i$,
for $ i = 1, \ldots , s$.  But then we obtain a homomorphism over $R$ 
using the Claim above. This proves the Corollary \ref{7c} and
hence the Proposition \ref{6p} in the case where $R$ is reduced and 
strictly henselian.  

\noindent
If $R$ is not reduced one applies standard deformation theory to $R
\rightarrow R_{\rm red}$,  [Z2], 4.47. 
\hfill $Q.E.D.$ 

\medskip

{\bf Proof} of Corollary \ref{7c}. Consider the diagram: 
\begin{equation*}
\begin{array}{ccc}
R & \rightarrow & R^{\rm sh}\\
\downarrow & & \downarrow\\
\hat{R} & \rightarrow & \hat{R}^{\rm sh}. 
\end{array}   
\end{equation*}
The upper index ``sh'' denotes the strict henselization. 
Using the fact that the categories of finite \'etale coverings of $R$,
$k$, respectively $\hat{R}$ are equivalent it is easy to see that the
last diagram is a fiber product. We have already proved that $p^c\psi$
lifts to a homomorphism over $R^{\rm sh}$. Applying Lemma \ref{7l} to the 
following diagram we see that $p^c\psi$ lifts to $\hat{R}$. This is
enough to prove the corollary (compare the Claim above). 
\begin{equation*}
\begin{array}{ccc}
R/\ideal{m}^n & \rightarrow & R^{\rm sh}/ \ideal{m}^n R^{\rm sh} \\
\downarrow & & \downarrow\\
R/ \ideal{m} & \rightarrow & R^{\rm sh}/ \ideal{m} R^{\rm sh} 
\end{array}   
\end{equation*}
In this diagram $n$ is a positive integer and $\ideal{m}$ is
the maximal ideal of $R$.  \hfill $Q.E.D.$
\vspace{-0.8cm}

\subsection*{}
\begin{thm} {\bf Corollary.} 
Let $R$ be a strictly henselian reduced local ring over
$\bar{\mathbb{F}}_p$. Let $R^{\rm perf}$ be the perfect hull of $R$. 
Let $X$ be a $p$-divisible group over
$S = \Spec R$ with constant Newton polygon. We set $S^{\rm perf} = 
\Spec(R^{\rm perf})$. Then there is a
$p$-divisible group $X_0$ over $\bar{\mathbb{F}}_p$ and an isogeny
$X_0 \times_{\Spec \bar{\mathbb{F}}_p} S^{\rm perf} \rightarrow X \times_S
S^{\rm perf}$ such that  the degree of this  isogeny is
bounded by an integer which depends only on the height of $X$.
\end{thm}
{\bf Proof:} This follows using Proposition \ref{perf}. 
\hfill $Q.E.D.$

\section{Examples}
\noindent
In this section we use the $p$-divisible groups either $Z = G_{1,n}$, 
with $n \geq 1$, or $Z = G_{m,1}$, with $m \geq 1$ as building blocks 
for our examples.  These have the property to be iso-simple, they are
defined over $\mathbb{F}_p$, they contain a unique subgroup scheme  $N \subset Z$
isomorphic with $\alpha_p$, and  $Z/N \cong Z$. 
Indeed, for  $Z = G_{1,n}$ we have an exact sequence of sheaves: 
$$ 0 \rightarrow \alpha_p \rightarrow Z \stackrel{\Fr}{\longrightarrow} Z
\rightarrow 0.$$
For $Z = G_{m,1}$ we have the exact sequence 
$$ 0 \rightarrow \alpha_p \rightarrow Z \stackrel{\Ver}{\longrightarrow} Z
\rightarrow 0.$$
Moreover every such $Z$ has the following property: if $Z_K \to Z'$ 
is an isogeny over some field $K$, and $k$  an algebraic 
closed field containing $K$, then $Z_k \cong Z'_k$.

\vspace{-0.8cm}

\subsection*{}
\begin{thm}\label{exa1} {\bf Example.}
In this example we produce 
a $p$-divisible group $X$ with constant Newton polygon over a regular
base scheme which does not admit a slope filtration.
\end{thm}
Choose $Z_1$ and $Z_2$ as above, with slope$(Z_1) = \lambda_1 > 
\lambda_2 = $ slope$(Z_2)$; 
e.g. $Z_1 = G_{1,1}$ and $Z_2=G_{1,2}$. We choose $R = K[t]$, 
where $K$ is a field. We write $S= \Spec(R)$ and  
$\mathcal{Z}_i = Z_i \times S$ for $i = 1, 2$. We define 
$$(id,t): \alpha_p \to   \alpha_p \times \alpha_p \cong N_1 
\times N_2; \quad\mbox{this defines}\quad \psi: \alpha_p
\times S \to \mathcal{Z}_1 \times \mathcal{Z}_2.$$

\smallskip
 
{\bf Claim:} $\mathcal{X} := (\mathcal{Z}_1 \times
\mathcal{Z}_2)/\psi(\alpha_p\times S)$ {\it is a $p$-divisible group
over $S$ which does not admit a slope filtration.}
 
\smallskip

Indeed, for the generic point we do have slope filtration, where 
$X =\mathcal{X} \otimes K(t)$, and $0 \subset X_1 \subset X$ is given by: 
$X_1$ is the image of 
$$\xi_K: (\mathcal{Z}_1 \otimes K(t) \to 
(\mathcal{Z}_1 \times \mathcal{Z}_2) \otimes K(t) \to X).$$ 
However the inclusion $\xi_K$ extends uniquely a homomorphism 
$\xi: \mathcal{Z}_1 \to \mathcal{X}$, which is not injective at $t \mapsto 0$. 
This proves the claim.        

\vspace{-0.8cm}

\subsection*{}
\begin{thm}\label{exa2} {\bf Example.} In this example we construct 
a $p$-divisible group $\mathcal{X}$ with constant Newton polygon over
a  base scheme $S$ \emph{which is  
not normal},  such that there is no isogeny $\phi: \mathcal{X} \to
\mathcal{Y}$ to a completely slope divisible $p$-divisible group.
(i.e. we show the condition that $S$ 
is  \emph{normal} in Theorem \ref{th} is necessary). 
\end{thm}
We start again with the exact sequences over $\mathbb{F}_p$:
\begin{equation}\label{exa2e}
   \begin{array}{cllll}
0 & \rightarrow \alpha_p & \rightarrow G_{2,1} &
\stackrel{\Ver}{\longrightarrow} G_{2,1} & \rightarrow 0,\\
0 & \rightarrow \alpha_p & \rightarrow G_{1,2} &
\stackrel{\Fr}{\longrightarrow} G_{1,2} & \rightarrow 0.\\
   \end{array}
\end{equation}
We fix an algebraically closed field $k$. We write $T = \mathbb{P}^1_k$, and:
\begin{displaymath}
   Z_1 = G_{2,1} \times_{\mathbb{F}_p} T, \quad 
   Z_2 = G_{1,2} \times_{\mathbb{F}_p} T, \quad 
   Z   = Z_1 \times Z_2, \quad 
   A   = \alpha_p \times_{\mathbb{F}_p} T.
\end{displaymath}
By base change we obtain sequences of sheaves on the projective line
$T = \mathbb{P}^1_k$:
\begin{displaymath}
   \begin{array}{cllll}
0 & \rightarrow A & \rightarrow Z_1 &
\stackrel{\Ver}{\longrightarrow} Z_1 & \rightarrow 0,\\
0 & \rightarrow A & \rightarrow Z_2 &
\stackrel{\Fr}{\longrightarrow} Z_2 & \rightarrow 0.\\
   \end{array}
\end{displaymath}

\vspace{3mm} \noindent
{\bf Lemma.}
{\it Consider $Z \to T = \mathbb{P}^1_k$ as above. Let $\beta: Z \rightarrow Y$ 
be an isogeny to a completely slope divisible
$p$-divisible group $Y$ over $T$. Then $Y = Y_1 \times Y_2$ is
a product of two $p$-divisible groups and $\beta = \beta_1 \times
\beta_2$ is the product of two isogenies $\beta_i : Z_i \rightarrow
Y_i$.} 
 
\vspace{3mm} \noindent
{\bf Proof.} The statement is clear if we replace the base $T$
by a perfect field, see Proposition \ref{perf}. In our case we show first
that the kernel of the morphisms $Z_i \rightarrow Y$, $i = 1,2$
induced by $\beta$ are representable by a finite, locally free group
schemes $G_i$. Indeed, let $\mathcal{G}_i$ the kernel in the sense of
f.p.p.f sheaves. Let us denote by $G$ the kernel of the isogeny
$\beta$. Choose a number $n$ such that $p^n$ annihilates $G$. Then
$p^n$ annihilates $\mathcal{G}_i$. Therefore $\mathcal{G}_i$ coincides
with the kernel of the morphism of finite group schemes $Z_i (n)
\rightarrow Y(n)$. Hence $\mathcal{G}_{i}$ is representable by a
finite group scheme $G_i$. We prove that $G_i$ is locally free. It
suffices to verify that the rank of $G_i$ in any geometric point $\eta$
of $T$ is the same. But we have seen that over $\eta$
the $p$-divisible group $Y$ splits into a product $Y_{\eta} =
(Y_{\eta})_1 \times (Y_{\eta})_2$. This implies that
\begin{displaymath}
   G_{\eta} = (G_1)_{\eta} \times (G_2)_{\eta}.
\end{displaymath}
We conclude that the ranks of $(G_i)_{\eta}$ are independent of $\eta$
since $G$ is locally free. 
 
\smallskip
 
We define $p$-divisible groups $Y_i = Z_i/G_i$. We obtain a
homomorphism of $p$-divisible groups $Y_1 \times Y_2 \rightarrow Y$
which is an isomorphism over each geometric point $\eta$. Therefore
this is an isomorphism.    \hfill $Q.E.D.$

\medskip

Next we construct a $p$-divisible group $X$ on $T = \mathbb{P}^1_k$. Let
$\mathcal{L}$ be a line bundle on $\mathbb{P}^1_k$. We consider the
associated vector group 
\begin{displaymath}
   \mathbf{V}(\mathcal{L}) (T') = \Gamma (T', \mathcal{L}_T'),
\end{displaymath}
where $T' \rightarrow T = \mathbb{P}^1_k$ is a scheme and
$\mathcal{L}_T'$ is the 
pull-back. The kernel of the Frobenius morphism $\Fr :
\mathbf{V}(\mathcal{L}) \rightarrow \mathbf{V}(\mathcal{L})^{(p)}$ is a
finite, locally free group scheme $\alpha_p (\mathcal{L})$ which is
locally isomorphic to $\alpha_p$. We set $A(-1) = \alpha_p
(\mathcal{O}_{\mathbb{P}^1_k}(-1))$. There are up to multiplication by
an element of $k^{\ast}$ unique homomorphisms 
$\iota_0$ respectively $\iota_{\infty} : \mathcal{O}_{\mathbb{P}^1_k}(-1)
\rightarrow \mathcal{O}_{\mathbb{P}^1_k}$ whose unique zeroes are $0 \in
\mathbb{P}^1_k$ respectively $\infty \in \mathbb{P}^1_k$. This induces
homomorphisms 
of finite group schemes $\iota_0 : A(-1) \rightarrow A$ respectively
$\iota_{\infty} : A(-1) \rightarrow A$ which are isomorphisms outside
$0$ respectively outside $\infty$. We consider the embeddings
\begin{displaymath}
   A(-1) \quad\stackrel{(\iota_0, \iota_{\infty})}{\longrightarrow}\quad 
   A \times A \quad\subset\quad Z_1 \times Z_2 = Z.
\end{displaymath}
We define $X = Z/A(-1)$:
$$\psi: Z \longrightarrow  Z/A(-1) = X.$$
Note that
$$\psi_0 : G_{2,1} \times G_{1,2} = Z_0 \longrightarrow G_{2,1} 
\times (G_{1,2}/\alpha_p) = X_0,$$
and
$$\psi_{\infty} : G_{2,1} \times G_{1,2} = Z_{\infty} 
\longrightarrow (G_{2,1}/\alpha_p) 
\times G_{1,2} = X_{\infty}.$$

\smallskip

We consider the quotient space $S = \mathbb{P}^1_k / \{ 0, \infty \}$, by 
identifying $0$ and $\infty$ into a normal crossing at $P \in S$, i.e. 
$$\mathcal{O}_{S,P} = \{f \in \mathcal{O}_{T,0} \cap
\mathcal{O}_{T,\infty} \mid f(0) = f(\infty)\};$$ 
$S$ is a nodal curve, and
$$ \mathbb{P}^1_k = T \longrightarrow S, \quad 0 \mapsto P, \quad
\infty \mapsto P,$$ 
is the normalization morphism.

\smallskip

A finite, locally
free scheme $\mathcal{G}$ over $S$ is the same thing as a finite, 
locally free scheme $G$ over $\mathbb{P}^1_k$ endowed with an isomorphism
\begin{displaymath}
   G_0 \cong G_{\infty}.
\end{displaymath}
It follows that the category of $p$-divisible groups $\mathcal{Y}$ over $S$ is
equivalent to the category of pairs $(Y, \gamma_{\mathcal{Y}})$, where $Y$ is
a $p$-divisible group on $\mathbb{P}^1_k$ and $\gamma_{\mathcal{Y}}$ is an
isomorphism 
$$\gamma_{\mathcal{Y}} : Y_0 \cong Y_{\infty}$$ of the fibers of $Y$
over $0 \in \mathbb{P}^1_k$ and $\infty \in \mathbb{P}^1_k$. We call
$\gamma_{\mathcal{Y}}$ the gluing datum of $\mathcal{Y}$. 
 
\smallskip
 
We construct  a  $p$-divisible $\mathcal{X}$ over $S$ by defining a gluing
datum on the $p$-divisible group $X$. In fact, the exact sequences 
in (\ref{exa2e}) give:
$$X_0 =  G_{2,1} 
\times (G_{1,2}/\alpha_p) \quad\cong\quad   (G_{2,1}/\alpha_p) 
\times G_{1,2}   = X_{\infty};  $$
this gluing datum provides a $p$-divisible group $\mathcal{X}$ over
$S$. 
 
\smallskip

{\bf Claim:} {\it This $p$-divisible group $\mathcal{X} \to S$ satisfies the 
property mentioned in the example.}

\smallskip

Let us assume that there exists an isogeny $\phi: \mathcal{X} \rightarrow
\mathcal{Y}$ to a completely slope divisible $p$-divisible group
$\mathcal{Y}$ over $S$.  We set $Y = \mathcal{Y} \times_{S} T$
and consider the induced isogeny 
$$ \phi \times_S T = \varphi:  X \rightarrow Y.$$
and the induced isogeny:
$$\beta = \varphi{\cdot}\psi: (Z \stackrel{\psi}{\longrightarrow} X 
\stackrel{\varphi}{\longrightarrow} Y).$$
By the lemma we see that
$$\beta = \beta_1 \times \beta_2 : Z \longrightarrow Y_1 
\times Y_2 = Y.$$
Note that 
$$\varphi_{\infty} = \phi_P = \varphi_0:   X_{\infty}=\mathcal{X}_P=X_0  \quad\longrightarrow 
\quad  Y_{\infty}=\mathcal{Y}_P=Y_0;$$
these $p$-divisible groups both have a splitting into isoclinic summands:
$$X_{\infty}=\mathcal{X}_P=X_0 = X' \times X'', \quad
Y_{\infty}=\mathcal{Y}_P=Y_0 = Y' \times  Y'',$$
and
$$\phi_P = \varphi' \times \varphi'': X' \times X''
\quad\longrightarrow\quad Y' \times  Y''$$
is in diagonal form.
On the one hand we conclude from
$$((\beta_1)_0: (Z_1)_0 \longrightarrow (Y_1)_0) = \left((Z_1)_0 
\stackrel{\sim}{\rightarrow}  X' \to Y'\right)$$
that $\deg(\beta_1) = \deg(\beta_1)_0 = \deg(\varphi');$
on the other hand 
$$((\beta_1)_{\infty}: (Z_1)_{\infty} \longrightarrow (Y_1)_{\infty}) = 
\left((Z_1)_{\infty} \to (G_{2,1}/\alpha_p)=X' \to Y'\right);$$
hence 
$$\deg(\beta_1) = \deg(\beta_1)_{\infty} = p{\cdot}\deg(\varphi').$$ 
We see that 
the assumption that the isogeny $\phi: \mathcal{X} \to \mathcal{Y}$ to a completely 
slope divisible $\mathcal{Y} \to S$ 
would exist leads to a contradiction. 
This finishes the description and the proof of  Example \ref{exa2}. 
\hfill $Q.E.D.$

\vspace{-0.8cm}

\subsection*{}
\begin{thm} {\bf Example.}  For every positive integer $d$
there exists a scheme $S'$ of dimension $d$, a point $P' \in S'$ such
that $S'$ is regular outside $P'$, and a $p$-divisible group
${\mathcal{X}}' \to S'$  which does not admit an isogeny to  a
completely slope divisible group over $S'$.
\end{thm}

This follows directly form the previous example. Indeed   choose $T$
as in the previous example, and let $T' \to T$ smooth and surjective
with $T'$ of dimension $d$.  Pull back  $X/T$  to  
$X'/T'$; choose geometric points $0'$ and $\infty' \in T'$ above $0$
and $\infty \in T$; construct $S'$ by ``identifying $0'$ and
$\infty'$'': outside  $P' \in S'$, 
this scheme is $T' \backslash \{0',  \infty'\}$, and the local ring of
$P' \in S'$ 
is the set of pairs of elements in the local rings of $0'$ and $\infty'$ having
the same residue value. We can descend $X' \to T'$
to ${\mathcal{X}}' \to S'$, and this has the desired property.

\Addresses
\end{document}